\documentclass[10pt]{article}
\usepackage{amsmath}

\usepackage{amsthm}
\theoremstyle{plain}
\newtheorem{thm}{Theorem}[section]

\newtheorem{lem}{Lemma}[section]

\theoremstyle{definition}
\newtheorem{defn}{Definition}[section]
\theoremstyle{plain}
\newtheorem{con}{conjecture}[section]
\newtheorem{rem}{Remark}[section]
\numberwithin{equation}{section}

\title{The maximal principle for properly immersed submanifolds and its applications}
\author{Yong Luo}
\date{}
\begin{document}
\maketitle

\begin{abstract}
In this note we consider the Liouville type theorem for a properly immersed submanifold $M$ in a complete Riemmanian manifold $N$. Assume that the sectional curvature $K^N$ of $N$ satisfies $K^N\geq-L(1+dist_N(\cdot,q_0)^2)^\frac{\alpha}{2}$ for some $L>0, 2>\alpha\geq 0$ and $q_0\in N$.

(i) If $\Delta|\vec{H}|^{2p-2}\geq k|\vec{H}|^{2p}$($p>1$) for some constant $k>0$, then we prove that $M$ is minimal.

(ii) Let $u$ be a smooth nonnegative function on $M$ satisfying $\Delta u\geq ku^a$ for some constant $k>0$ and $a>1$. If $|\vec{H}|\leq C(1+dist_N(\cdot,q_0)^2)^\frac{\beta}{2}$ for some $C>0$, $0\leq\beta<1$, then $u=0$ on $M$.

As applications we get some nonexistence result for $p$-biharmonic submanifolds.
\end{abstract}

\section{Introduction}

In the past several decades harmonic maps play a central role in geometry and analysis.
Let $\phi:(M^{m},g)\rightarrow(N^{m+t},h)$ be a map between Riemannian manifolds $(M, g)$ and $(N, h)$. The energy of $\phi$ is defined by
   \begin{equation*}
     E(\phi)=\int_{M}\frac{|d\phi|^{2}}{2}d\nu_{g},
   \end{equation*}
   where $d\nu_g$ is the volume element on $(M, g)$.

The Euler-Lagrange equation of $E$ is
\begin{equation*}
  \tau(\phi)=\sum_{i=1}^{m}\{ \tilde{\nabla}_{e_{i}}d\phi(e_{i})-d\phi(\nabla_{e_{i}}e_{i})\}=0,
\end{equation*}
where $\tilde{\nabla}$ is the Levi-Civita connection on the pullback bundle $\phi^{-1}TN$ and $\{e_i\}$ is a local orthonormal frame field on $M$.
In 1983, Eells and Lemaire \cite{EL} proposed to consider the bienergy functional
 \begin{equation*}
  E_{2}(\phi)=\int_{M}\frac{|\tau(\phi)|^{2}}{2}d\nu_{g},
\end{equation*}
where $\tau(\phi)$ is the tension field of $\phi$. Recall that $\phi$ is harmonic if $\tau(\phi)=0$. The Euler-Lagrange equation for $E_{2}$ is
  \begin{equation*}
    \tau_{2}(\phi)=\tilde{\triangle}(\tau(\phi))-\sum_{i=1}^{m}R^{N}(\tau(\phi),d\phi(e_{i}))d\phi(e_{i})=0.
  \end{equation*}
  To further generalize the notion of harmonic maps, Peter and Moser\cite{HM}(see also \cite{HF}) considered the $p(p>1)$-bienergy functional as follows:
\begin{equation*}
   E_{p}(\phi)=\int_{M}|\tau(\phi)|^{p}d\nu_{g}.
\end{equation*}
The $p$-bitension field $\tau_{p}(\phi)$ is
    \begin{equation}\label{bitensionfild}
\tau_{p}(\phi)=\tilde{\triangle}( |\tau(\phi)| ^{p-2}\tau(\phi))-\sum_{i=1}^{m}\bigg(R^{N}
  \big( |\tau(\phi)| ^{p-2}\tau(\phi),d\phi(e_{i})\big)d\phi(e_{i})\bigg).
    \end{equation}
The Euler-Lagrange equation for $E_p$ is $\tau_{p}(\phi)=0$ and a map $u$ satisfying $\tau_{p}(\phi)=0$ is called $p$-biharmonic maps. If $\phi:(M^{m},g)\rightarrow(N^{m+t},h)$ is an isometry immersion, then we call $u$ $p$-biharmonic submanifold and $2$-biharmonic submanifolds are called biharmonic submanifolds.

For biharmonic submanifolds, we have the well known Chen's conjecture \cite{Ch}:
\begin{con}
Every biharmonic submanifold in $\mathbf{E}^{n}$ is minimal.
\end{con}
Chen's conjecture inspires the research on the nonexistence of biharmonic submanifolds in nonpositively curved manifolds (\cite{Ak-Ma}\cite{BMO}\cite{CMO1} \cite{CMO2}\cite{CMO3}\cite{Ch} \cite{Ch1}\cite{CM} \cite{De}\cite{Di}\cite{Fu}\cite{Ha-Vl}\\\cite{Ji1}\cite{Luo1}\cite{Luo2}\cite{Ma}\cite{Na-Ur1}
\cite{Na-Ur2}\cite{NUG}  \cite{On} \cite{Ou} \cite{Ou-Ta} etc.). Motivated by Chen's conjecture, Yingbo Han \cite{Han} proposed the following conjecture:
\begin{con}\label{c2}
Every complete $p$-biharmonic submanifolds in nonpositively curved Riemannian manifold is minimal.
\end{con}
Some partial affirmative answers to conjecture \ref{c2} were proved in \cite{Han} and \cite{CL}. In this note we will continue to consider the nonexistence of  $p$-biharmonic submanifolds in nonpositively curved Riemannian manifold. Before mentioning our main result, we define the following notion(see \cite{Ma}).
\begin{defn}
For a complete manifold $(N, h)$ and $\alpha\geq0$, if the sectional curvature $K^N$ of $N$ satisfies $$K^N\geq -L(1+dist_N(\cdot, q_0)^2)^\frac{\alpha}{2},$$
for some $L>0$ and $q_0\in M$, then we say that $K^N$ has a polynomial growth bound of order $\alpha$ from below.
\end{defn}
We have
\begin{thm}\label{thm1}
Let $(M, g)$ be a properly immersed submanifold in a complete Riemannian manifold $(N, h)$ whose sectional curvature $K^N$ has a polynomial growth bound of order less than 2 from below. Assume that there exists a positive constant $k>0$ such that($p>1$)
\begin{eqnarray}
\Delta|\vec{H}|^{2p-2}\geq k|\vec{H}|^{2p}~ on ~M.
\end{eqnarray}
Then $M$ is minimal.
\end{thm}
\begin{rem}
When $p=2$, theorem \ref{thm1} was proved by Maeta(see \cite{Ma}). Our proof follows his argument by using the second derivatives' test to our new test functions. Maeta's argument was developed by Cheng and Yau in the 1970s(see \cite{CY1}\cite{CY2}\cite{CY3} etc.).
\end{rem}
Theorem \ref{thm1} implies the following nonexistence result of $p$-biharmonic submanifolds.
\begin{thm}\label{thm2}
Let $(M, g)$ be a properly immersed $p$-biharmonic submanifold in a complete nonpositively curved Riemannian manifold $(N, h)$ whose sectional curvature $K^N$ has a polynomial growth bound of order less than 2 from below, then $M$ is minimal.
\end{thm}
Using the same argument, we also have the following Liouville type theorem.
\begin{thm}\label{thm3}
Let $(M, g)$ be a properly immersed submanifold in a complete Riemannian manifold $(N, h)$ whose sectional curvature $K^N$ has a polynomial growth bound of order less than 2 from below. Assume that $u$ is a smooth nonnegative function on $M$ satisfying
\begin{eqnarray}
\Delta u\geq ku^a~on ~M,
\end{eqnarray}
where $k>0, a>1$ are constants. If $|\vec{H}|\leq C(1+dist_N(\cdot,q_0)^2)^\frac{\beta}{2}$ for some $C>0$, $0\leq\beta<1$ and $q_0\in N$, then $u=0$ on $M$. Here $\vec{H}$ is the mean curvature vector field of $M$ in $N$.
\end{thm}
This Liouville type theorem was first found by Maeta. In \cite{Ma} he proved the case of $a=2$.

 The rest of this paper is organized as follows: In section 2 we will briefly recall the theory of $p$-biharmonic submanifolds and submanifold theory. Our main theorems are proved in section 3.
\section{Preliminaries}

In this section we give more details on the definitions of harmonic maps, biharmonic maps, $p$-biharmonic maps and $p$-biharmonic
submanifolds.

Let $u:(M^{m},g)\rightarrow(N^{m+t},h)$ be a map from an $m$-dimensional Riemannian manifold $(M,g) $ to an $m+t$-dimensional Riemannian manifold $(N,h)$. The energy of $u$ is defined by
   \begin{equation*}
     E(u)=\int_{M}\frac{|du|^{2}}{2}d\nu_{g}.
   \end{equation*}
The Euler-Lagrange equation of $E$ is
\begin{equation*}
  \tau(u)=\sum_{i=1}^{m}\{ \tilde{\nabla}_{e_{i}}du(e_{i})-du(\nabla_{e_{i}}e_{i})\}=0,
\end{equation*}
where we denote $\nabla$ the Levi-Civita connection on $(M,g)$, and $\tilde{\nabla}$ the induced Levi-Civita connection of the pullback bundle $u^{-1}TN$. A map $u:(M^{m},g)\rightarrow(N^{m+t},h)$  is called a harmonic map if $\tau(u)=0$. To generalize the notion of harmonic maps, Eells and Lemaire \cite{EL} proposed to consider the bienergy functional
 \begin{equation*}
  E_{2}(u)=\int_{M}\frac{|\tau(u)|^{2}}{2}d\nu_{g}.
\end{equation*}
The Euler-Lagrange equation for $E_{2}$ is(see \cite{Ji2})
  \begin{equation*}
    \tau_{2}(u)=\tilde{\triangle}(\tau(u))-\sum_{i=1}^{m}R^{N}(\tau(u),du(e_{i}))du(e_{i})=0.
  \end{equation*}
  To further generalize the notion of harmoic maps, Han and Feng \cite{HF}(see also \cite{HM})introduced
the $F$-bienergy functional
\begin{equation*}
  E_{F}(u)=\int_{M}F(\frac{|\tau(u)|^{2}}{2})d\nu_{g},
\end{equation*}
where $F:[0, +\infty)$ and $F'(x)>0$ if $x>0$.

The critical points of the $F$-bienergy functional with $F(x)=(2x)^\frac{p}{2}(p>1)$ are called $p$-biharmonic maps and isometric $p$-biharmonic maps are called $p$-biharmonic submanifolds.

The $p$-bitension field $\tau_{p}(u)$ is
    \begin{equation}\label{bitensionfild}
      \tau_{p}(u)=\tilde{\triangle}( |\tau(u)| ^{p-2}\tau(u))-\sum_{i=1}^{m}\bigg(R^{N}
  \big( |\tau(u)| ^{p-2}\tau(u),du(e_{i})\big)du(e_{i})\bigg).
    \end{equation}
A $p$-biharmonic map satisfies $ \tau_{p}(u)=0$.

Now we briefly recall the submanifold theory. Let $u:(M,g)\rightarrow(N,h)$  be an isometric immersion from an m-dimensional Riemannian manifold into an $m + t$-dimensional Riemannian manifold. The second fundamental form $B: TM\times TM\to  T^\perp(M)$ is defined by:
 \begin{equation}
  \bar{\nabla}_{X}Y=\nabla_{X}Y+B(X,Y),X,Y\in \Gamma(TM),
 \end{equation}
 where $\bar{\nabla}$ is the Levi-Civita connection on $N$ and $\nabla$ is the Levi-Civita connection on $M$. The Weingarten formula is given by
\begin{equation}\label{connection  local and global}
   \bar{\nabla}_{X}\xi=-A_{\xi}X+\nabla^{\bot}_{X}\xi,X\in\Gamma(TM),,
\end{equation}
where $A_{\xi}$ is called the Weingarten map w.r.t. $\xi\in T^\perp(M)$, and $\nabla^{\bot}$ denotes the normal connection on the normal bundle of M in N. For any $x \in M$, the mean curvature vector field $\vec{H}$ of $M$ at $x$ is
 \begin{equation*}
   \vec{H}=\frac{1}{m}\sum_{i=1}^{m}B(e_{i},e_{i}).
 \end{equation*}
If $u$ is an isometry immersion, we see that $\{du(e_{i})\}$ is a local orthonormal frame of M. In addition, for any $X,Y\in \Gamma(TM)$,
 \begin{equation}
   \nabla du(X,Y)=\tilde{\nabla}_X(du(Y))-du(\nabla_{X}^{Y})=B(X, Y),
 \end{equation}
where $\tilde{\nabla}$ is the connection on the  pull back bundle $u^{-1}TN$, whose fiber at a point $x\in M$ is $T_{u(x)}N=T^{\bot}M\bigoplus TM$. Therefore if $u$ is an isometric immersion,
  \begin{equation*}
    \tau(u)=tr\nabla du=trB=m\vec{H},
  \end{equation*}
  and a $p$-biharmonic submanifold satisfies the following equatuion:
\begin{equation}\label{bitensionfild}
      \tau_{p}(u)=\tilde{\triangle}( |\vec{H}| ^{p-2}\vec{H})-\sum_{i=1}^{m}\bigg(R^{N}(
  \big |\vec{H}| ^{p-2}\vec{H}, e_{i}\big)e_{i}\bigg)
    \end{equation}
where $\tilde{\triangle}=\sum_{i=1}^{m}(\tilde{\nabla}_{e_{i}}\tilde{\nabla}_{e_{i}}-\tilde{\nabla}_{\nabla_{e_{i}e_{i}}})$, $\tilde{\nabla}$ is the connection on the pullback bundle,  and $R^{N}$ is the Riemanian curvature tensor on $N$.

From (\ref{connection  local and global}), we get for any vector field $\xi\in\Gamma(T^{\bot}M )$:
\begin{eqnarray*}
\tilde{\nabla}_{e_{i}}\tilde{\nabla}_{e_{i}}\xi
 &=&\tilde{\nabla}_{e_{i}}(\nabla^{\bot}_{e_{i}}\xi-A_{\xi}e_{i})\\
 &=&\nabla^{\bot}_{e_{i}}\nabla^{\bot}_{e_{i}}\xi-\tilde{\nabla}_{e_{i}}A_{\xi}e_{i}-A_{\nabla^{\bot}_{e_{i}}\xi}e_{i}\\
 &=&\nabla^{\bot}_{e_{i}}\nabla^{\bot}_{e_{i}}\xi-\nabla_{e_{i}}A_{\xi}e_{i}-B(e_{i},A_{\xi}e_{i})+A_{\nabla^{\bot}_{e_{i}}\xi}(e_{i}),
\end{eqnarray*}
and
 \begin{equation*}
   \begin{split}
    &\tilde{\nabla}_{\nabla_{e_{i}}e_{i}}\xi\\
    &=\nabla^{\bot}_{\nabla_{e_{i}}e_{i}}\xi-A_{\xi}(\nabla_{e_{i}}e_{i})
    \end{split}.
 \end{equation*}
Combining the above two identities, we get
  \begin{equation*}
  \begin{split}
   &\tilde{\triangle}\xi=\nabla^{\bot}_{e_{i}}\nabla^{\bot}_{e_{i}}\xi-\nabla_{e_{i}}A_{\xi}e_{i}-B(e_{i},A_{\xi}e_{i})+A_{\nabla^{\bot}_{e_{i}}\xi}(e_{i})\\
   &+\nabla^{\bot}_{\nabla_{e_{i}}e_{i}}\xi-A_{\xi}(\nabla_{e_{i}}e_{i})\\
   &=\triangle^{\bot}\xi-\nabla_{e_{i}}A_{\xi}e_{i}-A_{\xi}(\nabla_{e_{i}}e_{i})-B(e_{i},A_{\xi}e_{i})+A_{\nabla^{\bot}_{e_{i}}\xi}(e_{i})
  \end{split}
\end{equation*}
Therefore by decomposing the $p$-biharmonic submanifold equation into its normal and tangential parts respectively we get \cite{Han}:
  \begin{equation}\label{1}
  \Delta^{\bot} \bigg(|\vec{H}|^{p-2}\vec{H}\bigg)- \sum_{i=1}^{m}B(A_{|\vec{H}|^{p-2}\vec{H}}e_{i},e_{i})+\sum_{i=1}^{m}\big(R^{N}
  (|\vec{H}|^{p-2}\vec{H},e_{i})e_{i}\big)^{\bot}=0,
\end{equation}

\begin{equation}\label{2}
  Tr_{g}(\nabla A_{|\vec{H}|^{p-2}\vec{H}})+Tr_{g}[A_{\nabla^{\perp}|\vec{H}|^{p-2}\vec{H}}(.)]
  -\sum_{1}^{m}(R^{N}(|\vec{H}|^{p-2}\vec{H},e_{i})e_{i})^{\top}=0.
\end{equation}

\section{Proof of theorems}
In this section, we will need the following Hessian comparison theorem(see \cite{CLN}).
\begin{lem}\label{lem1}
Let $(N, h)$ be a complete Riemannian manifold with $sect\geq K(K<0)$. For any point $q\in M$ the distance function $r(x)=d(x,q)$ satisfies
$$D^2r\leq\sqrt{|K|}coth(\sqrt{|K|}r)h,$$
at all points where $r$ is smooth (i.e. away from $q$ and the cut loss). Here $D^2r$ denotes the Hessian of $r$.
\end{lem}
\subsection{Proof of Theorem \ref{thm1}}

\proof If $M$ is compact we see that $\vec{H}=0$ follows from the standard maximal principle. Therefore we assume that $M$ is noncompact. We will prove the theorem by a contradiction argument. Here we follow Maeta's(\cite{Ma}) argument by choosing new test functions.

Suppose that $\vec{H}(x_0)\neq 0$ for some $x_0\in M$. Set $u(x)=|\vec{H}(x)|^{2p-2}$ for $x\in M$. For each $\rho>0$ let
$$F(x)=F_\rho(x)=(\rho^2-r^2(\phi(x)))^{2p-2}u(x),$$
for $x\in M\cap X^{-1}(\bar{B}_\rho)$, where $\phi: M\to R^n$ is the isometric immersion, $B_\rho$ is the standard ball in $R^n$ with radius $\rho$ and $r(\phi(x))=dist_N(\phi(x), q_0)$ for some $q_0\in N$.

Assume that $x_0\in X^{-1}(B_{\rho_0})$. For each $\rho\geq\rho_0$, $F=F_\rho$ is a nonnegative function which is not identically zero on $M\cap X^{-1}(\bar{B}_\rho)$ and equals zero on the boundary. Assume that $q\in M\cap X^{-1}(B_\rho)$ is the maximum point of $F$($q$ exists because $\phi$ is properly immersed).

(i) $\phi(q)$ is not on the cut loss of $q_0$. Then $\nabla F(q)=0$ and hence we get at $q$
\begin{eqnarray}\label{F2}
\frac{\nabla u}{u}=\frac{(2p-2)\nabla r^2(\phi(x))}{\rho^2-r^2(\phi(x))}.
\end{eqnarray}
 In addition at $q$
 \begin{eqnarray}\label{F1}
 0\geq\Delta F(x)&=&(2p-2)(2p-3)(\rho^2-r^2(\phi(x)))^{2p-4}|\nabla r^2(\phi(x))|^2u(x)\nonumber
 \\&-&(2p-2)(\rho^2-r^2(\phi(x))^{2p-3}\Delta r^2(\phi(x))u(x)\nonumber
 \\&-&2(2p-2)(\rho^2-r^2(\phi(x))^{2p-3}\langle\nabla r^2(\phi(x)),\nabla u\rangle_g\nonumber
 \\&+&(\rho^2-r^2(\phi(x)))^{2p-2}\Delta u.
 \end{eqnarray}
 Combining  inequalities (\ref{F2}) and (\ref{F1}) we have at $q$
\begin{eqnarray}\label{ine1}
\frac{\Delta u(x)}{u(x)}\leq\frac{(2p-2)(2p-1)|\nabla r^2(\phi(x))|^2}{(\rho^2-r^2(\phi(x)))^2}
+\frac{(2p-2)\Delta r^2(\phi(x))}{\rho^2-r^2(\phi(x))}.
\end{eqnarray}
By a direct computation we see that
$$|\nabla r^2(\phi(x))|_g^2\leq 4m r^2(\phi(x)),$$
and
\begin{eqnarray}
&& \Delta r^2(\phi(x))
=2\sum_{i=1}^m\langle(\bar{\nabla}r)(\phi(x)), d\phi(e_i)\rangle^2\nonumber
\\&+&2r(\phi(x))\sum_{i=1}^m(D^2r)(\phi(x))\langle d\phi(e_i), d\phi(e_i)\rangle+2r(\phi(x))\langle(\bar{\nabla}r)(\phi(x)), \tau(\phi)(x)\rangle \nonumber
\\&\leq&2m+2r(\phi(x))\sum_{i=1}^m(D^2r)(\phi(x))\langle d\phi(e_i), d\phi(e_i)\rangle+2mr(\phi(x))|\vec{H}(x)|,
\end{eqnarray}
 where $m=dimM, \bar{\nabla}$ is the gradient on $(N, h)$ and $D^2r$ denotes the Hessian of $r$. Since the sectional curvature $K^N$ of $N$ satisfies $K^N\geq-L(1+r^2)^{\frac{\alpha}{2}}$, by the Hessian comparison theorem(see lemma \ref{lem1}) we get
 \begin{eqnarray}
 \sum_{i=1}^m(D^2r)(\phi(x))\langle d\phi(e_i), d\phi(e_i)\rangle\leq m\sqrt{L(1+r^2)^{\frac{\alpha}{2}}}coth\bigg(\sqrt{L(1+r^2)^{\frac{\alpha}{2}}}r(\phi(x))\bigg).
 \end{eqnarray}
 Combining the last two inequalities we obtain
 \begin{eqnarray}\label{ine2}
 \Delta r^2(\phi(x))&\leq &2m+2m\sqrt{L(1+r^2)^{\frac{\alpha}{2}}}r(\phi(x))coth\bigg(\sqrt{L(1+r^2)^{\frac{\alpha}{2}}}r(\phi(x))\bigg)\nonumber
 \\&+&2mr(\phi(x))|\vec{H}(x)|.
 \end{eqnarray}Recall that $\Delta|\vec{H}|^{2p-2}\geq k|\vec{H}|^{2p}$, i.e. $\Delta u\geq ku^\frac{2p}{2p-2},$ thus from inequalities (\ref{ine1}), (\ref{ine2}) we obtain
\begin{eqnarray}
&&ku(q)^\frac{1}{p-1}\leq \frac{4m(2p-2)(2p-1)r^2(\phi(q))}{(\rho^2-r^2(\phi(q)))^2}\nonumber
\\&+&\frac{(2p-2)\bigg\{2m
+2m\sqrt{L(1+r^2)^{\frac{\alpha}{2}}}r(\phi(q))coth\bigg(\sqrt{L(1+r^2)^{\frac{\alpha}{2}}}r(\phi(q))\bigg)\bigg\}}{\rho^2-r^2(\phi(q))}\nonumber
\\&+&\frac{(2p-2)2mr(\phi(q))|\vec{H}(q)|}{\rho^2-r^2(\phi(q))}.
\end{eqnarray}
From the last inequality one gets
\begin{eqnarray}
&&u(q)\leq C(p,k,m)[\frac{r^{2p-2}(\phi(q))}{(\rho^2-r^2(\phi(q)))^{2p-2}}\nonumber
\\&+&\frac{\bigg\{1
+\sqrt{L(1+r^2)^{\frac{\alpha}{2}}}r(\phi(q))coth\bigg(\sqrt{L(1+r^2)^{\frac{\alpha}{2}}}r(\phi(q))\bigg)\bigg\}^{p-1}}{(\rho^2-r^2(\phi(q)))^{p-1}}\nonumber
\\&+&\sqrt{u(q)}r(\phi(q))^{p-1}\frac{1}{(\rho^2-r^2(\phi(q)))^{p-1}}],
\end{eqnarray}
where $C(p,k,m)$ is a constant depends only on $p,k,m$. Therefore
\begin{eqnarray}
&&F(q)\leq C(p,k,m)[r^{2p-2}(\phi(q))+\nonumber
\\&&\bigg\{1
+\sqrt{L(1+r^2)^{\frac{\alpha}{2}}}r(\phi(q))coth\bigg(\sqrt{L(1+r^2)^{\frac{\alpha}{2}}}r(\phi(q))\bigg)\bigg\}^{p-1}(\rho^2-r^2(\phi(q)))^{p-1} \nonumber
\\&+&\sqrt{F(q)}r(\phi(q))^{p-1}]
\end{eqnarray}
which implies that
$$F(q)\leq C(p,k,m,L)(1+\rho^2)^{\frac{(\alpha+6)}{4}(p-1)},$$
where $C(p,k,m,L)$ is a constant depends only on $p,k,m,L$.
Since $q$ is the maximum of $F$, for any $x\in M\cap B_\rho$ we have
$$F(x)\leq F(q)\leq C(p,k,m,L)(1+\rho^2)^{\frac{(\alpha+6)}{4}(p-1)}.$$
Therefore
\begin{eqnarray}
|\vec{H}(x)|^{2p-2}\leq \frac{C(p,k,m,L)(1+\rho^2)^{\frac{(\alpha+6)}{4}(p-1)}}{(\rho^2-r^2(\phi(x)))^{2p-2}},
\end{eqnarray}
for any $x\in M\cap B_\rho$ and $\rho\geq \rho_0$.

(ii) If $\phi(q)$ is on the cut loss of $q_0$, then we use a method of Calabi (see \cite{Ca}). Let $\sigma$ be a minimal geodesic joining $\phi(q)$ and $q_0$. Then for any $q'$ in the interior of $\sigma$, $q'$ is not conjugate to $q_0$. Fix for such a point $q'$. Let $U_{q'}\subseteq B_\rho$ be a conical neighborhood of the geodesic segment of $\sigma$ joining $q'$ and $\phi(q)$ such that for any $\phi(x)\in U_{q'}$, there is at most one minimizing geodesic joining $q'$ and $\phi(x)$. Let $\bar{r}(\phi(x))=dist_{U_{q'}}(\phi(x), q')$ in the manifold $U_{q'}$. Then we have $\bar{r}(\phi(x))\geq dist_N(\phi(x),q')$, $r(\phi(x))\leq r(q')+\bar{r}(\phi(x))$, $r(\phi(q))=r(q')+\bar{r}(\phi(q))$. We claim that the function
$$F_{\rho,q'}(x):=(\rho^2-\{r(q')+\bar{r}(\phi(x))\}^2)^{2p-2}u(x)~for ~x\in \phi^{-1}(U_{q'})$$
also attains a local maximum at the point $q$. In fact, for any point $x\in \phi^{-1}(U_{q'})$ we have
\begin{eqnarray*}
F_{\rho,q'}(q)&=&(\rho^2-\{r(q')+\bar{r}(\phi(q))\}^2)^{2p-2}u(q)
\\&=&(\rho^2-r^2(\phi(q)))^{2p-2}u(q)
\\&=&F_\rho(q)\geq F_\rho(x)
\\&=&(\rho^2-r^2(\phi(x)))^{2p-2}u(x)
\\&\geq &(\rho^2-\{r(q')+\bar{r}(\phi(x))\}^2)^{2p-2}u(x)
\\&=&F_{\rho,q'}(x).
\end{eqnarray*}
Therefore the claim is proved and we play the second derivative's test to $F_{\rho,q'}(x)$ at $q$, the same argument as before shows that
\begin{eqnarray*}
&&F_{\rho,q'}(q)\leq C(p,k,m)[r^{2p-2}(\phi(q))+\nonumber
\\&&\bigg\{1
+\sqrt{L(1+r^2)^{\frac{\alpha}{2}}}r(\phi(q))coth\bigg(\sqrt{L(1+r^2)^{\frac{\alpha}{2}}}r(\phi(q))\bigg)\bigg\}^{p-1}(\rho^2-r^2(\phi(q)))^{p-1} \nonumber
\\&+&\sqrt{F_{\rho,q'}(q)}r(\phi(q))^{p-1}],
\end{eqnarray*}
which implies that
$$F_{\rho,q'}(q)\leq C(p,k,m,L)(1+\rho^2)^{\frac{(\alpha+6)}{4}(p-1)}.$$
Take $q'\to q_0$ we have $F_{\rho,q'}(q)=F_\rho(q)$ and hence
$$F_\rho(q)\leq C(p,k,m,L)(1+\rho^2)^{\frac{(\alpha+6)}{4}(p-1)}.$$
Therefore
\begin{eqnarray}
|\vec{H}(x)|^{2p-2}\leq \frac{C(p,k,m)(1+\rho^2)^{\frac{(\alpha+6)}{4}(p-1)}}{(\rho^2-r^2(\phi(x)))^{2p-2}},
\end{eqnarray}
for any $x\in M\cap B_\rho$ and $\rho\geq \rho_0$. Let  $x=x_0$ and $\rho\to +\infty$ we get $\vec{H}(x_0)=0,$  a contradiction. Therefore $M$ is minimal.
\endproof
\subsection{Proof of Theorem \ref{thm2}}
\proof Recall that the normal part of the $p$-biharmonic submanifolds is
  \begin{eqnarray*}\label{Nor}
  \Delta^{\bot} \bigg(|\vec{H}|^{p-2}\vec{H}\bigg)- \sum_{i=1}^{m}B(A_{|\vec{H}|^{p-2}\vec{H}}e_{i},e_{i})+\sum_{i=1}^{m}\big(R^{N}
  (|\vec{H}|^{p-2}\vec{H},e_{i})e_{i}\big)^{\bot}=0.
\end{eqnarray*}
Therefore
\begin{eqnarray*}
\Delta|\vec{H}|^{2p-2}&=&2\langle\Delta^\bot(|\vec{H}|^{p-2}\vec{H}),|\vec{H}|^{p-2}\vec{H}\rangle+2|\nabla(|\vec{H}|^{p-2}\vec{H})|^2
\\&\geq&2\sum_{i=1}^{m}\langle B(A_{|\vec{H}|^{p-2}\vec{H}}e_{i},e_{i}), |\vec{H}|^{p-2}\vec{H}\rangle
\\&=&2|\vec{H}|^{2p-4}\langle A_{\vec{H}}e_i, A_{\vec{H}}e_i\rangle
\\&\geq& 2m|\vec{H}|^{2p},
\end{eqnarray*}
where in the first inequality we used the assumption of nonpositive curvature. Therefore $M$ is minimal by theorem \ref{thm1}.
\endproof
\subsection{Proof of Theorem \ref{thm3}}
\proof Similar to the proof of theorem \ref{thm1} set $F_\rho(x)=(\rho^2-r^2(\phi(x)))^{2a-2}u(x)$. If $u(x_0)\neq0$, then using the second derivatives' test to $F_\rho$ at the maximum point $q$ for $\rho$ big enough such that $x_0\in B_\rho$, we will get
\begin{eqnarray}
&&ku(q)^\frac{1}{a-1}\leq \frac{4m(2a-2)(2a-1)r^2(\phi(q))}{(\rho^2-r^2(\phi(q)))^2}\nonumber
\\&+&\frac{(2a-2)\bigg\{2m
+2m\sqrt{L(1+r^2)^{\frac{\alpha}{2}}}r(\phi(q))coth\bigg(\sqrt{L(1+r^2)^{\frac{\alpha}{2}}}r(\phi(q))\bigg)\bigg\}}{\rho^2-r^2(\phi(q))}\nonumber
\\&+&\frac{(2a-2)2mr(\phi(q))|\vec{H}(q)|}{\rho^2-r^2(\phi(q))}.
\end{eqnarray}
Therefore \begin{eqnarray}
&u(q)\leq C(a,k,m)[\frac{r^{2a-2}(\phi(q))}{(\rho^2-r^2(\phi(q)))^{2a-2}}\nonumber
\\&+\frac{\bigg\{1
+\sqrt{L(1+r^2)^{\frac{\alpha}{2}}}r(\phi(q))coth\bigg(\sqrt{L(1+r^2)^{\frac{\alpha}{2}}}r(\phi(q))\bigg)\bigg\}^{a-1}}{(\rho^2-r^2(\phi(q)))^{a-1}}\nonumber
\\&+|\vec{H}|^{a-1}r(\phi(q))^{a-1}\frac{1}{(\rho^2-r^2(\phi(q)))^{a-1}}],
\end{eqnarray}
which implies that
$$F_\rho(q)\leq  C(a,k,m,L)\max\{(1+\rho^2)^\frac{\alpha+6}{4}, (1+\rho^2)^\frac{\beta+3}{2}\},$$
where we used the assumption that $|\vec{H}|\leq C(1+dist_N(\cdot,q_0)^2)^\frac{\beta}{2}$.
Therefore
$$(\rho^2-r^2(\phi(x)))^{2a-2}u(x)\leq F_\rho(q)\leq C(a,k,m,L)\max\{(1+\rho^2)^{\frac{\alpha+6}{4}(a-1)}, (1+\rho^2)^{\frac{\beta+3}{2}(a-1)}\},$$
which implies that
$$u(x)\leq\frac{C(a,k,m,L)\max\{(1+\rho^2)^{\frac{\alpha+6}{4}(a-1)}, (1+\rho^2)^{\frac{\beta+3}{2}(a-1)}\}}{(\rho^2-r^2(\phi(x)))^{2a-2}}.$$
Because $\alpha<2$ and $\beta<1$, let $x=x_0$ and $\rho\to +\infty$ we obtain $u(x_0)=0$, a contradiction. Thus $u=0$ on $M$.
\endproof
\textbf{Acknowledgement:} The author would like to thank Prof. Yingbo Han for stimulating discussions and Dr. Shun Maeta for his kind suggestions and comments.

\vspace{1cm}\sc
  Yong Luo

School of mathematics and statistics,

Wuhan university, Hubei 430072, China

{\tt yongluo@whu.edu.cn}

\vspace{1cm}\sc
\end{document}